\documentclass{amsart}

\usepackage{graphicx}
\usepackage{wrapfig}
\usepackage{lscape}
\usepackage{rotating}
\usepackage{epstopdf}
\usepackage{comment}
\usepackage{amssymb}

\usepackage{enumerate}
\usepackage{enumitem}
\usepackage{longtable}
\usepackage{geometry}
\usepackage{tabularx}

\usepackage{subcaption}

\usepackage{mdframed}
\usepackage{xcolor}
\usepackage{graphicx}
\usepackage{amsmath}
\usepackage{bm}
\usepackage{amsthm}
\usepackage{thmtools}
\usepackage{float}
\usepackage{breqn}

\newtheorem{theorem}{Theorem}[section]
\newtheorem{lemma}[theorem]{Lemma}
\newtheorem{conjecture}[theorem]{Conjecture}
\newtheorem{proposition}[theorem]{Proposition}
\newtheorem{corollary}[theorem]{Corollary}

\theoremstyle{definition}
\newtheorem{definition}[theorem]{Definition}
\newtheorem{example}[theorem]{Example}

\theoremstyle{remark}
\newtheorem{remark}[theorem]{Remark}

\newcommand{\cor}{\begin{corollary}}
	\newcommand{\roc}{\end{corollary}}
\newcommand{\pr}{\begin{proof}}
	\newcommand{\rp}{\end{proof}}
\newcommand{\ex}{\begin{example}}
	\newcommand{\xe}{\end{example}}
\newcommand{\tm}{\begin{theorem}}
	\newcommand{\mt}{\end{theorem}}
\newcommand{\lm}{\begin{lemma}}
	\newcommand{\ml}{\end{lemma}}
\newcommand{\df}{\begin{definition}}
	\newcommand{\fd}{\end{definition}}
\newcommand{\prop}{\begin{proposition}}
	\newcommand{\porp}{\end{proposition}}
\newcommand{\cj}{\begin{conjecture}}
	\newcommand{\jc}{\end{conjecture}}
\newcommand{\rmk}{\begin{remark}}
	\newcommand{\kmr}{\end{remark}}

\theoremstyle{plain}
\newtheorem*{GGC1}{Generalised Goldbach conjecture (GGC)}

\theoremstyle{plain}

\numberwithin{equation}{section}

\def\modd#1 #2{#1\ \mbox{\rm (mod}\ #2\mbox{\rm )}}

\begin{document}

\title[]{Predicting the size ranking of minimal primes in the generalised Goldbach partitions}

\author{Zs\'ofia Juh\'asz}
\address{Dept. of Computer Algebra \\
	Faculty of Informatics\\
	E\"otv\"os Lor\'and University \\
	1117 Budapest \\
	Hungary}
\curraddr{}
\email{jzsofia@inf.elte.hu}
\thanks{}

\author{M\'at\'e Bartalos}
\address{}
\curraddr{}
\email{bartmate@gmail.com}
\thanks{}


\date{}

\dedicatory{To my Grandmothers}

\begin{abstract}
	 A scarcely known  generalization of Goldbach's conjecture  introduced by Hardy and Littlewood  states that for every pair of (relatively prime) positive integers $m_1$ and $m_2$, every sufficiently large integer $n$ satisfying certain simple congruence criteria can be $(m_1, m_2)$-partitioned as $n=m_1p+m_2q$ for some primes $p$ and $q$.  While the size of the minimal prime in the Goldbach partitions of even numbers has received prior attention, we extend this investigation to the general case of $(m_1, m_2)$-partitions.  This question has a direct implication on the running times of verification algorithms of the generalised Goldbach conjecture.  We study the rankings of the pairs $(m_1, m_2)$   according to the sizes of the averages and maxima, respectively, of the minimal $p$ in the $(m_1, m_2)$-partitions of numbers up to large thresholds, and propose a rank-order predicting function depending only on  $m_2$ and the prime factors of $m_1$. We computed both the average and the maximum of the minimal prime $p$ in all $(m_1, m_2)$-partitions of integers up to $10^9$, for every pair of relatively prime coefficients  $1\leq m_1 \neq m_2\leq 40$. Our function shows very high rank-order correlations with both the empirical averages and maxima of the minimal primes $p$ (Spearman’s $\rho =0.9949$ and $0.9958$, respectively). It also correctly predicts trends in the experimental data, for example, that for all relatively prime  $1\leq m_1 < m_2\leq 40$, the average minimal $p$ in the $(m_1, m_2)$-partitions of numbers up to $10^9$ exceeds the analogous average for the $(m_2, m_1)$-partitions.	 We present numerical data, including the average and the maximum of the minimal  $p$ in the $(m_1, m_2)$-partitions of numbers up to $10^9$ for each pair  $1\leq m_1 \neq m_2\leq 20$ relatively prime, and the resulting size rankings.
\end{abstract}

\maketitle

\section{Introduction}	\label{Sect_Intro}

The even (or strong) Goldbach  conjecture poses one of the longest-standing open problems in number theory,  which has been intriguing wide audiences for nearly three centuries. Goldbach formulated his famous claim  in 1742 in a letter to Euler, proposing  --  rephrased in modern form --   that every even number greater than $2$ can be expressed as the sum of two primes \cite{Fuss}.

Although a general proof has been elusive, major relevant results have been achieved. In 1973 Chen proved  that every sufficiently large even number is the sum of a prime and a semiprime (the product of at most two primes) \cite{Chen}. In 2010 Lu showed  that the number of even integers up to $x$ which do not have Goldbach partitions is $O(x^{0.879})$ \cite{Wen}.  In 2013 Helfgott proved the odd (weak or ternary) Goldbach conjecture  -- a related but weaker statement than the even Goldbach conjecture -- asserting that every odd number greater than 5 is the sum of three primes \cite{Helfgott2013, Helfgott2015}.

Several efforts have been made to verify the even Goldbach conjecture empirically up to increasing limits \cite{Pipping, Richstein, Sinisalo, Granville1989}. The current record is $4\cdot 10^{18}$ set by Oliveira e Silva \textit{et al.}\  in 2014 \cite{Oliveira2014}.

In 1894 Lemoine \cite{Lemoine} suggested a stronger version of the weak Goldbach conjecture, stating that every odd number $n > 5$ can be expressed as $n = p + 2q$ for some
primes $p$ and $q$  \cite[p.\ 424]{Dickson1919}. Although still unproven,  Lemoine's conjecture (LC) has drawn more limited research interest than those of  Goldbach. The highest limit of its empirical verification  was $10^{10}$ \cite{MaketheBrain} before 2024, when it was lifted to $10^{13}$ by the authors \cite{JuhaszBartalos2024}.

In 1923 Hardy and Littlewood \cite{Hardy} introduced the following generalization (GGC) of the even Goldbach conjecture, also generalizing LC:   for every (relatively prime) positive integer $m_1$ and $m_2$, every sufficiently large integer $n$ satisfying certain simple congruence conditions can be expressed as $n=m_1p+m_2q$ for some primes $p$ and $q$. Interestingly, this generalization appears to be hardly known and not to have been studied until  2017  \cite{FarkasJuhasz}. 
The authors recently tested GGC empirically up to $10^{12}$ [and up to $10^{13}$ for selected pairs] for all relatively prime $1\leq m_1,m_2\leq 40$ \cite{FarkasJuhasz, JuhaszBartalos2024}, presenting the largest values of $n$  found satisfying the conditions of GGC without $(m_1, m_2)$-partitions. The small sizes of these largest counterexamples support the plausibility of GGC. Lemoine's conjecture was also confirmed up to a new record of $10^{13}$.

The smallest value of $p$ in all Goldbach partitions of an even number $n\geq 4$ is denoted by $p(n)$. The rate of growth of $p(n)$ has been the subject of prior research [5] and directly influences the time complexity of algorithms verifying Goldbach’s conjecture by searching for the minimal partitions [5, 13]. Granville \textit{et al.}\  \cite{Granville1989}  conjectured   $p(n)=O(\log ^2 n\log \log n)$. Granville also suggested two more precise, incompatible conjectures of the form $p(n)\leq (C+o(1))\log^2 n\log{\log{n}}$, where $C$ is `sharp' in the sense that $C$ is the smallest constant with this property: one with $C=C_2^{-1}\approx 1.51478$ and the other one with $C=2e^{-\gamma}C_2^{-1}\approx 1.70098$, where $C_2\approx 0.66016$ is the twin prime constant and $\gamma\approx 0.57722$ is the Euler's constant \cite{Oliveira2014}. Empirical comparison of the plausibility of these conjectures was inconclusive due to the need for data up to even higher limits \cite{Oliveira2014}.

In the general setting of GGC,  for any integers $m_1, m_2$ and $n$ satisfying the conditions of GGC, we call an expression $n=m_1p+m_2q$ where $p$ and $q$ are primes, an $(m_1, m_2)$-partition of $n$.  Let $p^*_{m_1, m_2}(n)$ denote the smallest value of $p$  in all $(m_1, m_2)$-partitions of $n$ if such partition exists \cite{JuhaszBartalos2024}. The magnitude of  $p^*_{m_1, m_2}(n)$ affects the runtime of corresponding verification algorithms \cite{JuhaszBartalos2024}.
To the authors' knowledge  the only published results regarding $p^*_{m_1, m_2}(n)$  in general, beyond the classical case $p^*_{1, 1}(n)=p^*(n)$, are from 2024, when the average and maximum of the minimal  $p^*_{m_1, m_2}(n)$ up to $10^9$ were presented for every $1\leq m_1 \neq m_2\leq 20$ relatively prime  \cite{JuhaszBartalos2024}.

In this paper we investigate how the relative sizes of the (average) values of the function $p^*_{m_1, m_2}(n)$ (where $n$ is up to a sufficiently large threshold) depend on $m_1$ and $m_2$. We study how the ranking of pairs $(m_1, m_2)$ according to the size of the average values of $p^*_{m_1, m_2}(n)$ 
can be predicted based on $m_1$ and $m_2$ only. Guided by a heuristic argument, we introduce a function $R_{m_1, m_2}$ depending only on  $m_2$ and the prime factors of  $m_1$. This function exhibits excellent rank-order correlations with both the average and the maximum of $p^*_{m_1, m_2}(n)$ up to $10^9$ over the pairs of relatively prime $1\leq m_1\neq m_2\leq 40$. The corresponding Spearman’s rank-order coefficients are $0.9949$ ($4$ d.p.) for the averages and $0.9958$ ($4$ d.p.) for the maxima.  
The function $R_{m_1, m_2}$ also forecasts two notable observations emerging from our data: (1) For every pair $1\leq m_1 < m_2\leq 40$ relatively prime,  the average value of $p^*_{m_1, m_2}(n)$ where $n\leq 10^9$ is larger than that of $p^*_{m_2, m_1}(n)$; and (2) The average value  of $p^*_{m_1, m_2}(n)$ where $n\leq 10^9$, apart from $m_2$, appears to be  affected only by the distinct prime factors of $m_1$, but hardly at all  by the exact value of $m_1$.

After preliminaries,  the main results are presented in Section \ref{Sect_Main}. Section \ref{Sect_Time} includes remarks regarding the connection between the magnitude of  $p^*_{m_1, m_2}(n)$ and the running times of verification algorithms for GGC. Section \ref{Sect_Concl} summarizes some of the main conclusions. Section \ref{Sect_ToD} contains the tables of relevant numerical data generated, including the following: the average and maximum values  of $p^*_{m_1, m_2}(n)$ where $n\leq 10^9$, for every $1\leq m_1 \neq m_2\leq 20$ relatively prime, and the rankings of all pairs $1\leq m_1 \neq m_2\leq 20$ relatively prime according to the size of the average $p^*_{m_1, m_2}(n)$ and that of the maximal $p^*_{m_1, m_2}(n)$, respectively, where $n\leq 10^9$, and according to the value $R_{m_1, m_2}$.

\section{Preliminaries} \label{SectPrelim}
For every integer $a$ and $b$, let $\gcd{(a, b)}$ denote the greatest common divisor of $a$ and $b$. Hardy and Littlewood \cite{Hardy} introduced the following conjecture:
\begin{GGC1} 	\label{conRGGC} 
	Let $m_1$ and $m_2$ be positive integers such that $\gcd(m_1, m_2)=1$. Then for every sufficiently large integer $n$ satisfying the conditions
	\begin{enumerate}
		\item $\gcd(n, m_1)=\gcd(n, m_2)=1$ and   \label{conRGGC1}
		\item $n\equiv \modd{m_1+m_2} {2}$,   \label{conRGGC2}
	\end{enumerate}
	there exist primes $p$ and $q$ such that 
	\begin{equation} \label{EqnGGpartition}
		n=m_1p+m_2q.
	\end{equation}
\end{GGC1}
Furthermore, they also conjectured the following estimate for the number of ways $N(n)$ in which an integer $n$ satisfying the conditions of  GGC can be expressed in the form \ref{EqnGGpartition}:
$$ N(n)\sim \frac{2C_2}{m_1m_2}\frac{n}{(\log n)^2}\Pi \left(\frac{p-1}{p-2}\right),$$
where $C_2$ is the twin prime constant, and the product is taken over all odd primes $p$ which divide $m_1, m_2$ or $n$.

In the sequel  $n$, $m_1$, and $m_2$ denote positive integers such that $m_1$ and $m_2$ are relatively prime. The definitions and notation below were introduced earlier \cite{JuhaszBartalos2024}:

\begin{definition}
An expression of the form \ref{EqnGGpartition} where $p$ and $q$ are primes is called an $(m_1, m_2)$\textit{-Goldbach partition} (or $(m_1, m_2)$\textit{-partition}) \textit{of} $n$.   We say that 
$n$ \textit{can be} $(m_1, m_2)$\textit{-partitioned} if it possesses at least one $(m_1, m_2)$-partition.
\end{definition}
\begin{definition}
If  $n$ can be $(m_1, m_2)$-partitioned then the smallest and the largest values of $p$ [$q$] in all $(m_1, m_2)$-partitions of $n$ are denoted by $p^*_{m_1, m_2}(n)$ [$q^*_{m_1, m_2}(n)$] and $p^{**}_{m_1, m_2}(n)$  [$q^{**}_{m_1, m_2}(n)$], respectively. We call $n=m_1p^*_{m_1, m_2}(n)+m_2q^{**}_{m_1, m_2}(n)$ the $p$\textit{-minimal}  (or $q$\textit{-maximal}) and 
$n=m_1p^{**}_{m_1, m_2}(n)+m_2q^{*}_{m_1, m_2}(n)$ the $p$\textit{-maximal}  (or $q$\textit{-minimal}) $(m_1, m_2)$\textit{-partition of} $n$. 
\end{definition}

 We let GGC$_{m_1, m_2}$ denote the claim of GGC for given coefficients $m_1, m_2$. Then GGC$_{1, 1}$ and GGC$_{1, 2}$ are Goldbach's and Lemoine's conjectures, respectively.  
Clearly, for every $m_1$, $m_2$  the conditions  of GGC$_{m_1, m_2}$ and GGC$_{m_2, m_1}$ on $n$ are equivalent, and every $(m_1, m_2)$-partition of $n$ is also  an $(m_2, m_1)$-partition if the order of terms is disregarded. Hence, a number $n$  can be $(m_1, m_2)$-partitioned if and only if it can be $(m_2, m_1)$-partitioned, and in this case $p^*_{m_1, m_2}(n)=q^*_{m_2, m_1}(n)$ and $p^{**}_{m_1, m_2}(n)=q^{**}_{m_2, m_1}(n)$. Conjectures GGC$_{m_1, m_2}$ and GGC$_{m_2, m_1}$ are equivalent.

For every $m_1, m_2$, the number $n=m_1+m_2$  satisfies the conditions of GGC$_{m_1, m_2}$ and cannot be $(m_1, m_2)$-partitioned. Hence, if GGC$_{m_1, m_2}$ is true then there exists a largest positive integer satisfying the conditions of GGC$_{m_1, m_2}$  that  cannot be $(m_1, m_2)$-partitioned, which we denote by $k_{m_1, m_2}$. For every $1\leq m_1, m_2\leq 40$ relatively prime, the authors computed the  largest integer less than or equal to $ 10^{12}$ satisfying the conditions of GGC$_{m_1, m_2}$ that cannot be $(m_1, m_2)$-partitioned, and denoted it by $\hat{k}_{m_1, m_2}$ and hypothesised that $\hat{k}_{m_1, m_2}=k_{m_1, m_2}$. \cite{JuhaszBartalos2024}. Clearly, we have $\hat{k}_{m_1, m_2}=\hat{k}_{m_2, m_1}$, and if GGC$_{m_1, m_2}$ (and so  GGC$_{m_2, m_1}$) holds, then $k_{m_1, m_2}=k_{m_2, m_1}$.

For  every pair $m_1, m_2$ and every $L> \hat{k}_{m_1, m_2}$ for which there is at least one  $n$ satisfying the conditions of GGC$_{m_1, m_2}$ such that $\hat{k}_{m_1, m_2} <n\leq L$,  we refer to the average [maximum] value of $p^*_{m_1, m_2}(n)$ over all $\hat{k}_{m_1, m_2}< n\leq L$ satisfying the conditions of GGC$_{m_1, m_2}$ more succinctly as the \textit{average} [\textit{maximum}] \textit{of} $p^*_{m_1, m_2}(n)$ \textit{up to} $L$. For brevity, the average  $p^*_{m_1, m_2}(n)$ up to $L$ is also written as $\overline{p^*}_{m_1, m_2}$ \textit{up to} $L$.
Table \ref{Table_AvgMax} in Section \ref{Sect_ToD} shows the average and  maximum $p^*_{m_1, m_2}(n)$ up to $10^9$ for each pair of relatively prime coefficients $1\leq m_1 \neq m_2\leq 20$.

\section{Function $R_{m_1, m_2}$: Comparing the average and the maximum $p^*_{m_1, m_2}(n)$ values, respectively, for different pairs $m_1, m_2$} \label{Sect_Main}

In this section we introduce the function $R_{m_1, m_2}$ of $m_1$ and $m_2$.  
For every positive integer $n$,  the value $\varphi(n)$ of Euler's totient function at $n$ is the number of positive integers less than or equal to $n$ that are relatively prime to $n$. For every integer $a$ and $m\neq 0$,   the modulo $m$ residue of $a$ is denoted by $a \bmod m$, the modulo $m$ residue class of $a$ is denoted by $[a]_m$. The ring of modulo $m$ residue classes and the multiplicative group of the reduced modulo $m$ residue classes are denoted  by $\mathbb{Z}/m\mathbb{Z}$ and by $\mathbb{Z}'_m$, respectively. For every $[a]_m, [b]_m\in \mathbb{Z}/m\mathbb{Z}$, $[a]_m+ [b]_m$ and $[a]_m[b]_m$ respectively, stand for the sum and the product  of $[a]_m$ and $[b]_m$ in $\mathbb{Z}/m\mathbb{Z}$. If $(a, m)=1$ then $[a]^{-1}_m$ denotes the inverse of $[a]_m\in \mathbb{Z}_m^'$ in $\mathbb{Z}'_m$. The set of all primes is  denoted by  $P$.

Our heuristics is based on the following proposition:

\prop \label{Prop_qexistsforp}
Let $m_1, m_2$ and $n$ be pairwise relatively prime positive integers satisfying the conditions  of GGC. If $m_2\neq 1$ then denote the canonical decomposition of $m_2$ by $m_2=r_1^{\alpha_1}\cdot\ldots \cdot r_k^{\alpha_k}$. For every prime $p<\frac{n}{m_1}$ and $q$ the following two statements are equivalent:

\bigskip

\begin{enumerate}

 \item We have $n=m_1p+m_2q$.

\bigskip
\item One of the following two cases holds:

\smallskip

\noindent Case 1: The following conditions are all satisfied:

\smallskip
\begin{enumerate}[label=(\alph*)] \label{Prop_qexistsforp_1}
	\item $q\;\;\nmid \;\; m_2$ and $q\;\;\nmid \;\; m_1$; \label{Prop_qexistsforp_1_1}
	\item if $m_2\neq 1$ then for every $1\leq i\leq k$: $p$ belongs to one of the following $r_i-1$ different (reduced) modulo $r_i^{\alpha _i+1}$ residue classes:  $[n]_{r_i^{\alpha _i+1}}[m_1]^{-1}_{r_i^{\alpha _i+1}} + j_i\cdot [r_i^{\alpha_i}]_{r_i^{\alpha _i+1}}$, where  $1\leq j_i\leq  r_i-1$; \label{Prop_qexistsforp_1_2}
	\item $p$ belongs to one of the following $q-1$ different modulo $q^{2}$ residue classes:  $[n]_{q^2}\cdot [m_1]_{q^2}^{-1}+\ell\cdot [q]_{q^2}$, where $1\leq \ell \leq q-1$ and  \label{Prop_qexistsforp_1_3}
	\item for every prime $s<K$   -- where $K$  can be chosen as any number $K\geq \frac{n}{m_2}$ -- such that $s\neq q$, $s\;\;\nmid \;\; m_1$ and $s\;\;\nmid \;\; m_2$: $p$ belongs to one of the following $s-1$ different modulo $s$ residue classes: $[n]_{s}\cdot [m_1]_{s}^{-1}+[j]_s$, where $1\leq j\leq s-1$. \label{Prop_qexistsforp_1_4}
\end{enumerate}

\bigskip

\noindent Case 2: The following conditions are all  satisfied:

\smallskip

\begin{enumerate}[label=(\alph*),start=5] \label{Prop_qexistsforp_2}
	\item $q | m_2$ (i.e. $q=r_{\ell}$ for some $1\leq \ell \leq k$) and $q\;\;\nmid \;\; m_1$; \label{Prop_qexistsforp_2_1} 
	\item  for every $1\leq i\leq k$, $i\neq \ell$: $p$ belongs to one of the following $r_i-1$ different (reduced) modulo $r_i^{\alpha _i+1}$ residue classes:  $[n]_{r_i^{\alpha _i+1}}\cdot [m_1]_{r_i^{\alpha _i+1}}^{-1}+j_{i}\cdot [r_i^{\alpha _i}]_{r_i^{\alpha _i+1}}$, where  $1\leq j_i\leq  r_i-1$; \label{Prop_qexistsforp_2_2}
	\item $p$ belongs to one of the following $r_{\ell}-1$ different (reduced) modulo $r_{\ell}^{\alpha _{\ell}+2}$ residue classes: $[n]_{r_{\ell}^{\alpha _{\ell}+2}}\cdot [m_1]_{r_{\ell}^{\alpha _{\ell}+2}}^{-1}+j_{\ell}\cdot [r_{\ell}^{\alpha _{\ell}+1}]$, where $1\leq j_{\ell}\leq r_{\ell}-1$ and   \label{Prop_qexistsforp_2_3}
	\item for every prime $s<K$  -- where $K$  can be chosen as any number $K\geq \frac{n}{m_2}$ --  such that $s\;\nmid \; m_1$ and $s\;\nmid \; m_2$: $p$ belongs to one of the following $s-1$ different modulo $s$ residue classes: $[n]_{s}\cdot [m_1]_{s}^{-1}+[j]_s$, where $1\leq j\leq s-1$. \label{Prop_qexistsforp_2_4}
\end{enumerate}

\end{enumerate}

\porp

\bigskip

\pr
If $n=m_1p+m_2q$ then $q\;\nmid \; m_1$, since otherwise $q|m_1p$ and $q|m_2q$ would imply $q|n$, and so $q|\gcd{(n, m_1)}=1$, a contradiction. Therefore we have the following two possible cases: either ($q\;\nmid \; m_1$ and $q\;\nmid \; m_2$) or ($q\;\nmid \; m_1$ and $q|m_2$). Clearly, $n=m_1p+m_2q$ $\Leftrightarrow$  $n-m_1p=m_2q$.

Case I: Suppose $q\;\nmid \; m_1$ and $q\;\nmid \; m_2$. Then if $m_2=1$ then $m_2q=q$. If $m_2\neq 1$ then the canonical decomposition of $m_2q$ is $m_2q=qr_1^{\alpha_1}\cdot\ldots \cdot r_k^{\alpha_k}$. Since $0<n-m_1p<n$, we have $n-m_1p=m_2q$ if and only if the following conditions are satisfied:

\begin{enumerate}[label=(\roman*)] \label{PrCase1}
	\item  \label{PrCase1/1} if $m_2\neq 1$ then for every $1\leq i\leq k$: $r_i^{\alpha_i}|(n-m_1p)$ and $r_i^{\alpha_i+1} \;\nmid \; (n-m_1p)$,
	
	\item $q|(n-m_1p)$ and  $q^2 \;\nmid \; (n-m_1p)$  and \label{PrCase1/2}
	
	\item for every prime $s<\frac{n}{m_2}$ such that $s\neq q$ and $s\;\nmid \; m_2$: $s\;\nmid \; (n-m_1p)$.	\label{PrCase1/3}	
\end{enumerate}

We show that Conditions \ref{PrCase1/1}, \ref{PrCase1/2} and \ref{PrCase1/3}  above are equivalent to
Conditions \ref{Prop_qexistsforp_1_2}, \ref{Prop_qexistsforp_1_3} and \ref{Prop_qexistsforp_1_4}, respectively, of the Proposition.

\medskip

\ref{PrCase1/1}
For every $1\leq i\leq k$:

\begin{itemize}
	\item $r_i^{\alpha _i}|(n-m_1p)$ 
	$\Leftrightarrow$ $[n]_{r_i^{\alpha _i}}=[m_1]_{r_i^{\alpha _i}}[p]_{r_i^{\alpha _i}}$ $\Leftrightarrow$ $[p]_{r_i^{\alpha _i}}=[n]_{r_i^{\alpha _i}}[m_1]^{-1}_{r_i^{\alpha _i}}$.
	
	\item $r_i^{\alpha _i+1}\;\nmid \;  (n-m_1p)$ 
	$\Leftrightarrow$ $[n]_{r_i^{\alpha _i+1}} \neq [m_1]_{r_i^{\alpha _i+1}}[p]_{r_i^{\alpha _i+1}}$ $\Leftrightarrow$ $[p]_{r_i^{\alpha _i+1}} \neq [n]_{r_i^{\alpha _i+1}}[m_1]^{-1}_{r_i^{\alpha _i+1}}$.
\end{itemize}

As $([p]_{r_i^{\alpha _i}}=[n]_{r_i^{\alpha _i}}[m_1]^{-1}_{r_i^{\alpha _i}}  \textrm{ and }   [p]_{r_i^{\alpha _i+1}} \neq [n]_{r_i^{\alpha _i+1}}[m_1]^{-1}_{r_i^{\alpha _i+1}})$  $\Leftrightarrow$
$[p]_{r_i^{\alpha _i+1}} = [n]_{r_i^{\alpha _i+1}}[m_1]^{-1}_{r_i^{\alpha _i+1}} + j_i\cdot [r_i^{\alpha_i}]_{r_i^{\alpha _i+1}}$ for some $1\leq j_i\leq r_i-1$, Conditions \ref{PrCase1/1} and \ref{Prop_qexistsforp_1_2} are equivalent.

\medskip

\ref{PrCase1/2} We have

\begin{itemize}
	\item $q|(n-m_1p)$ 
	 $\Leftrightarrow$ $[m_1]_q[p]_q=[n]_q$ $\Leftrightarrow$ $[p]_q=[n]_q[m_1]^{-1}_q$.

	\item $q^2 \;\nmid \; (n-m_1p)$ 
	$\Leftrightarrow$ $[m_1]_{q^2}[p]_{q^2} \neq [n]_{q^2}$ $\Leftrightarrow$ $[p]_{q^2} \neq [n]_{q^2}[m_1]^{-1}_{q^2}$.
\end{itemize}

Since

$([p]_q=[n]_q[m_1]^{-1}_q \textrm{ and } [p]_{q^2} \neq [n]_{q^2}[m_1]^{-1}_{q^2})$ $\Leftrightarrow$ $[p]_{q^2}=[n]_{q^2}[m_1]^{-1}_{q^2}+\ell \cdot [q]_{q^2}$ for some $1\leq \ell \leq q-1$, Conditions \ref{PrCase1/2} and \ref{Prop_qexistsforp_1_3} are equivalent.

\medskip

\ref{PrCase1/3} 
If $s|m_1$ then -- since $\gcd{(n, m_1)}=1$ -- $s\;\nmid \;  n$, and so $s\;\nmid \;  (n-m_1p)$. Therefore Condition \ref{PrCase1/3} is equivalent to the following: for every prime $s$ such that $s\neq q$,  $s\;\nmid \;  m_1$ and $s\;\nmid \;  m_2$: $s\;\nmid \;  (n-m_1p)$.	For every prime $s$ such that $s\neq q$,  $s\;\nmid \;  m_1$ and $s\;\nmid \;  m_2$:

$s\;\nmid \;  (n-m_1p)$ $\Leftrightarrow$ $[n]_s \neq [m_1]_s[p]_s$ $\Leftrightarrow$ $[p]_s \neq [n]_s  [m_1]_s^{-1}$ $\Leftrightarrow$  $[p]_s = [n]_s  [m_1]_s^{-1}+ [j]_s$ for some $1\leq j \leq s-1$. Therefore Conditions \ref{PrCase1/3} and \ref{Prop_qexistsforp_1_4} are equivalent. Hence in Case I, $n=m_1p+m_2q$ is equivalent to Conditions \ref{Prop_qexistsforp_1_1}, \ref{Prop_qexistsforp_1_2}, \ref{Prop_qexistsforp_1_3} and \ref{Prop_qexistsforp_1_4} of the Proposition being all satisfied.

\medskip

Case II: Suppose $q\;\nmid \;  m_1$ and $q| m_2$. Then $q=r_{\ell}$ for some $1\leq \ell \leq k$ and the canonical decomposition of $m_2q$ is $m_2q=r_1^{\alpha_1}\cdot\ldots \cdot r_{\ell -1}^{\alpha_{\ell -1}}\cdot q^{\alpha_{\ell}+1} \cdot r_{\ell +1}^{\alpha_{\ell +1}} \cdot \ldots\cdot r_k^{\alpha_k}$.

Similarly to Case I it can be shown that in Case II, $n=m_1p+m_2q$  is equivalent to Conditions \ref{PrCase2/1}, \ref{PrCase2/2} and \ref{PrCase2/3} below being all satisfied, where Conditions \ref{PrCase2/1}, \ref{PrCase2/2} and \ref{PrCase2/3} are equivalent to Conditions \ref{Prop_qexistsforp_2_2}, \ref{Prop_qexistsforp_2_3} and \ref{Prop_qexistsforp_2_4}, respectively, of the Proposition:

\begin{enumerate}[label=(\roman*)]\setcounter{enumi}{3} \label{PrCase2}	
	\item for every $1\leq i\leq k$, $i\neq \ell$: $r_i^{\alpha_i}|(n-m_1p)$ and $r_i^{\alpha_i+1}\;\nmid \; (n-m_1p)$, \label{PrCase2/1}	
	\item $q^{\alpha_{\ell}+1}|(n-m_1p)$ and  $q^{\alpha_{\ell}+2}\;\nmid \; (n-m_1p)$  and \label{PrCase2/2}	
	\item for every prime $s<\frac{n}{m_2}$ such that $s\;\nmid \;  m_2$: $s\;\nmid \;  (n-m_1p)$.	\label{PrCase2/3}	
\end{enumerate}

 Hence in Case II, $n=m_1p+m_2q$ is equivalent to Conditions \ref{Prop_qexistsforp_2_1}, \ref{Prop_qexistsforp_2_2}, \ref{Prop_qexistsforp_2_3} and \ref{Prop_qexistsforp_2_4} of the Proposition being all satisfied.

\rp

Next we describe a heuristic argument to deduce the function $R_{m_1, m_2}$ of ordered pairs $(m_1, m_2)$ of relatively prime coefficients. In the following let $m_1, m_2$ and $n$ be given positive integers satisfying the conditions of GGC.

Condition \ref{Prop_qexistsforp_1_2}  in Case 1 of Proposition \ref{Prop_qexistsforp}  means that if $m_2\neq 1$ then $\forall 1\leq i\leq k$  the prime $p$ belongs to one of certain $r_i-1$ reduced modulo $r_i^{\alpha_i+1}$ classes. Hence, using the Chinese Remainder Theorem and the Prime Number Theorem for Arithmetic Progressions (PNTAP), for a given prime $q$, the ratio of primes satisfying this condition among all primes in a 'large' interval can be estimated as $$\prod_{1\leq i\leq k}\frac{r_i-1}{\varphi{(r_i^{\alpha_i+1})}}=\prod_{1\leq i\leq k}\frac{r_i-1}{r_i^{\alpha_i+1}-r_i^{\alpha_i}}=\prod_{1\leq i\leq k}\frac{1}{r_i^{\alpha_i}}=\frac{1}{m_2}.$$

 If $q|n$ then Condition \ref{Prop_qexistsforp_1_3} means that $q|(n-m_1p)$ and $q^2 \;\nmid \; n-m_1p$. Then $q|m_1p$, and since $n$ and $m_1$ are relatively prime, $q|p$ and so $q=p$. Therefore, in this case there is at most one $p$, namely $p=q$ which can satisfy Condition \ref{Prop_qexistsforp_1_3}. Hence, the ratio of those primes $p$ in a `large' interval satisfying the condition  tends  to $0$, as the length of the interval approaches infinity. If $q\;\nmid \; n$ then Condition \ref{Prop_qexistsforp_1_3} means that the prime $p$ belongs to one of certain $q-1$ reduced modulo $q^2$ residue classes, so in this case,  by PNTAP the ratio of  primes $p$ satisfying Condition \ref{Prop_qexistsforp_1_3} among all primes in a `large' interval can be estimated by $\frac{q-1}{\varphi(q^2)}=\frac{q-1}{q^2-q}=\frac{1}{q}$.

Let $q$ be a given prime. Condition \ref{Prop_qexistsforp_1_4} means that for every prime $s$  such that $s\neq q$, $s\;\nmid \; m_1$ and $s\;\nmid \; m_2$ we have the following: If $s|n$ then $p\;\nmid \;  s$, that is $p\neq s$, and if $s\;\nmid \;  n$ then either $p=s$ or $p$ belongs to one of certain $s-2$ reduced modulo $s$ residue classes.  The ratio of  primes $p$ satisfying Condition \ref{Prop_qexistsforp_1_4} among all primes in a `large' interval can be estimated by
$$\prod_{\substack{s\in P, s\neq q, s<\frac{n}{m_2} \\ s\;\nmid \; nm_1m_2  }}\frac{s-2}{\varphi(s)}=\prod_{\substack{s\in P, s\neq q,  s<\frac{n}{m_2} \\ s\;\nmid \; nm_1m_2 }}\frac{s-2}{s-1},$$ 
(Note that since $n\equiv m_1+m_2 \pmod{2}$, exactly one of $n, m_1$ and $m_2$ is even, and so $s\;\nmid \; nm_1m_2$ implies that in the factors of the above product $s\neq 2$.)

If for some primes $p$ and $q$ we have $q|n$ and $n=m_1p+m_2q$ then $q|m_1p$, which -- since $\gcd{(n, m_1)}=1$-- implies $q|p$, hence  $p=q$, and so $n=(m_1+m_2)q$. 
Therefore, if $q|n$, then ignoring the case of the single number $n=(m_1+m_2)q$, the number of primes satisfying  $n=m_1p+m_2q$ is $0$. Furthermore, since for given $m_1, m_2$, $n$ and $K$ no prime $p$ can satisfy Conditions \ref{Prop_qexistsforp_1_3} and \ref{Prop_qexistsforp_1_4}  simultaneously for more than one prime $q\leq K$, the ratio of those primes $p$ for which there exists a prime $q\leq \frac{n}{m_2}$, $q\;\nmid \;  m_1m_2$ such that Conditions \ref{Prop_qexistsforp_1_2} \ref{Prop_qexistsforp_1_3} and \ref{Prop_qexistsforp_1_4} are simultaneously satisfied, in a `large' interval can be estimated as
\begin{dmath} \label{Eqn1}
	\sum\limits_{\substack{q\in P, q< \frac{n}{m_2} \\ q\;\nmid \; n m_1m_2}} \left(\frac{1}{m_2q}\cdot \prod\limits_{\substack{s\in P, s\neq q \\ s\;\nmid \; nm_1m_2 \\ s<\frac{n}{m_2}}}\frac{s-2}{s-1}\right) = \frac{1}{m_2} \cdot \left( \prod_{\substack{s\in P, s<\frac{n}{m_2}\\  s\;\nmid \; nm_1m_2 }}\frac{s-2}{s-1} \right)\cdot \sum_{\substack{q\in P, q< \frac{n}{m_2} \\ q\;\nmid \;  nm_1m_2}}\frac{q-1}{q(q-2)}\approx\frac{1}{m_2} \cdot \left( \prod_{\substack{s\in P, s<\frac{n}{m_2}\\  s\;\nmid \; nm_1m_2 }}\frac{s-2}{s-1} \right)\cdot \sum_{\substack{q\in P, q< \frac{n}{m_2} \\ q\;\nmid \;  nm_1m_2}}\frac{1}{q-1}.
\end{dmath}

Case 2:  Estimated ratio of primes $p< \frac{n}{m_1}$ satisfying the conditions of Case 2 for some prime $q|m_2$, among all primes less than $\frac{n}{m_1}$: 
\begin{dmath}
		\sum\limits_{\substack{q\in P, q|m_2}} \left(\frac{1}{m_2q} \cdot \prod\limits_{\substack{s\in P, s\leq\frac{n}{m_2} \\ s\;\nmid \; nm_1m_2}}\frac{s-2}{s-1}\right) = \frac{1}{m_2} \cdot \left(\prod\limits_{\substack{s\in P, s\leq\frac{n}{m_2} \\ s\;\nmid \; nm_1m_2}}\frac{s-2}{s-1} \right) \cdot \sum\limits_{\substack{q\in P,  q|m_2}}\frac{1}{q} \approx\frac{1}{m_2} \cdot \left(\prod\limits_{\substack{s\in P, s\leq\frac{n}{m_2} \\ s\;\nmid \; nm_1m_2}}\frac{s-2}{s-1}\right) \cdot \sum\limits_{\substack{q\in P, q|m_2}}\frac{1}{q-1}.
\end{dmath}
The sum of the above estimates yields an estimate for the ratio of primes $p< \frac{n}{m_1}$ satisfying the conditions of Case 1 or Case 2 for some prime $q$ among all primes less than $\frac{n}{m_1}$:

\begin{equation}
\frac{1}{m_2} \cdot \left(\prod\limits_{\substack{s\in P, s\leq\frac{n}{m_2} \\ s \;\;\nmid \; \; nm_1m_2}}\frac{s-2}{s-1}\right) \cdot \sum\limits_{\substack{q\in P, q\leq\frac{n}{m_2}\\ q \;\;\nmid \; \; nm_1}} \frac{1}{q-1}
\end{equation}
For ease of calculation define
\begin{equation}
R'_{m_1, m_2}(n):=\frac{1}{m_2} \cdot \left(\prod\limits_{\substack{s\in P, s\leq n \\ s\;\nmid \;nm_1m_2}}\frac{s-2}{s-1}\right) \cdot \sum\limits_{\substack{q\in P, q\leq n\\ q\;\nmid \;nm_1}} \frac{1}{q-1}. 
\end{equation}
Then for all sufficiently large $n$ satisfying the conditions of GGC$_{m_1, m_2}$:
\begin{equation}
	\begin{split}
\frac{ R'_{1, 1}(n)}{ R'_{m_1, m_2}(n)}  & =   \left(\prod_{\substack{s\in P, s\;\nmid \; n \\ 2<s\leq n}}\frac{s-2}{s-1}\right) \cdot \left(\sum\limits_{\substack{q\in P, s\;\nmid \; n \\ q\leq n}}\frac{1}{q-1}\right) \bigg/\left(\frac{1}{m_2} \cdot \left(\prod_{\substack{s\in P, s\leq n \\ s\;\nmid \; nm_1m_2}}\frac{s-2}{s-1}\right) \cdot \left(\sum\limits_{\substack{q\in P, q\leq n \\ q\;\nmid \; nm_1}}\frac{1}{q-1}\right)  \right) \\ & =   m_2 \cdot \left(\prod_{\substack{s\in P, 2<s\leq n \\ s | m_1m_2}}\frac{s-2}{s-1}\right)\cdot\left(1+\left(\sum\limits_{\substack{q\in P, q\leq n \\  q | m_1}}\frac{1}{q-1}\right)/\left(\sum\limits_{\substack{q\in P, q\leq n \\ q \;\nmid \; nm_1}}\frac{1}{q-1}\right)  \right) 
\end{split}
\end{equation}

Define  $\displaystyle R_{m_1, m_2} := m_2 \cdot \prod_{\substack{s\in P, 2<s\\ s | m_1m_2}}\frac{s-2}{s-1}$.

We hypothesize that there is a `strong' rank-order correlation between the following three values across all pairs of relatively prime positive integers $(m_1, m_2)$: (1)  $\overline{p^*}_{m_1, m_2}$ up to any sufficiently large threshold;  
(2) the maximum of $p^*_{m_1, m_2}$ up to any sufficiently large threshold and  (3) $R_{m_1, m_2}$.

\begin{figure}[H]
	\centering
	\begin{subfigure}{0.325\textwidth}
		\includegraphics[width=\textwidth]{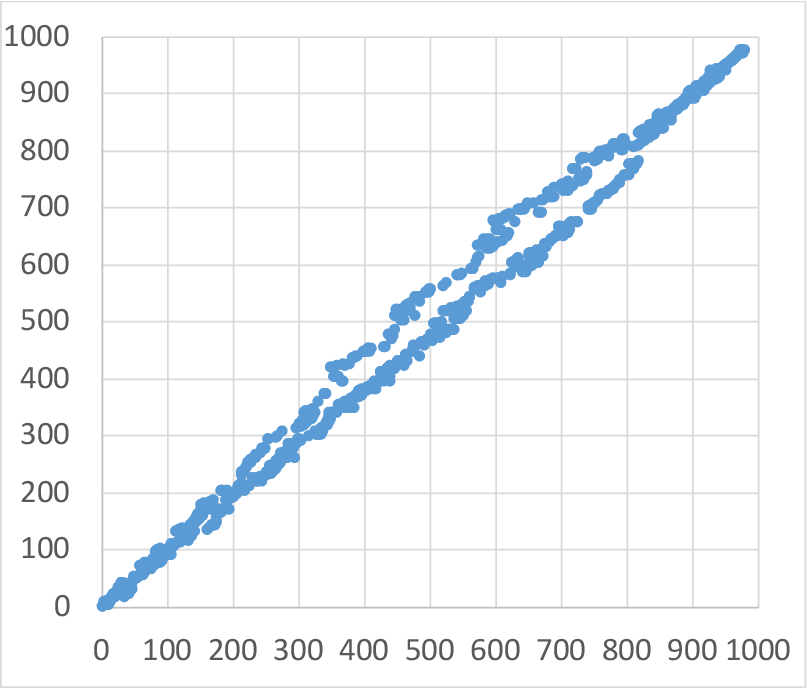}
		\caption{Ranking 3 (vertical ax.) plotted against Ranking 1}
		\label{subfig1}
	\end{subfigure}
	\hfill
	\begin{subfigure}{0.325\textwidth}
		\includegraphics[width=\textwidth]{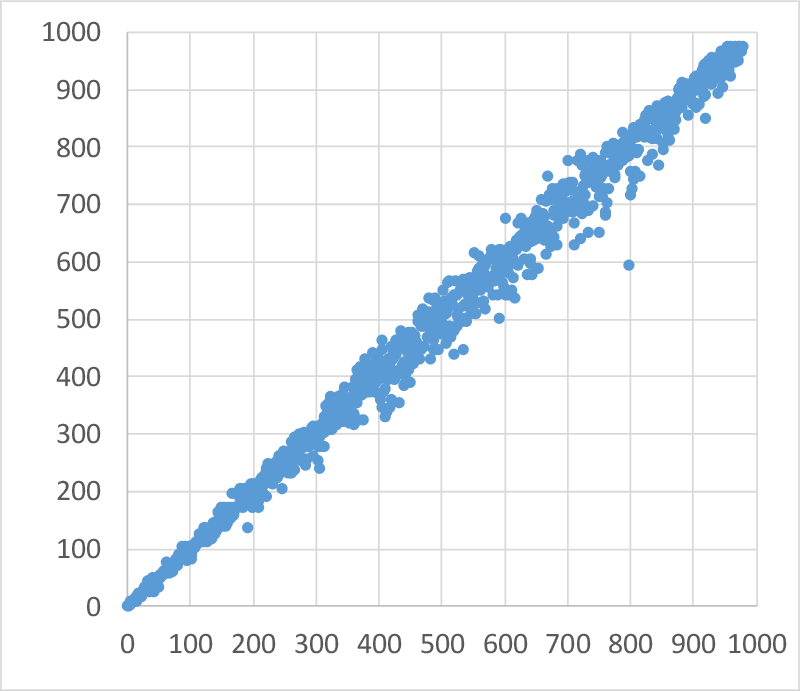}
		\caption{Ranking 3 (vertical ax.) plotted against Ranking 2}
		\label{subfig2}
	\end{subfigure}
	\hfill
	\begin{subfigure}{0.325\textwidth}
		\includegraphics[width=\textwidth]{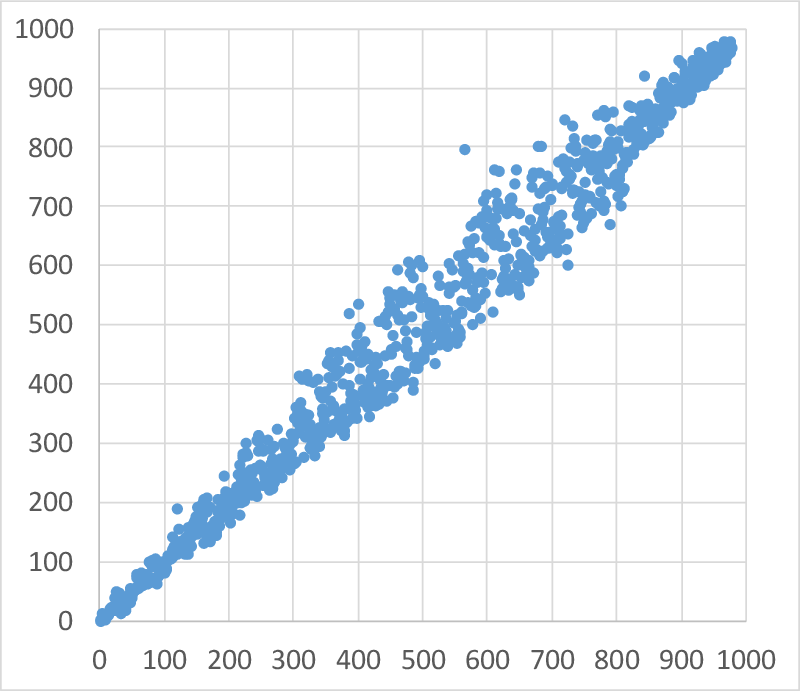}
		\caption{Ranking 2 (vertical ax.) plotted against Ranking 1}
		\label{subfig3}
	\end{subfigure}	
	\caption{Subfigures \ref{subfig1} and  \ref{subfig2} plot the ascending ranking of all pairs of relatively prime $1\leq  m_1 \neq m_2 \leq 40$ by $R_{m_1, m_2}$ on the vertical axis against their ascending ranking  by  $\overline{p^*}_{m_1, m_2}$  up to $10^9$ 	and against their ascending ranking by the maximum of $p^*_{m_1, m_2}$ up to $10^9$,   respectively. Subfigure \ref{subfig3} displays the ascending ranking of the same pairs ${m_1, m_2}$ by the maximum of $p^*_{m_1, m_2}$  up to $10^9$  against their ascending ranking by  $\overline{p^*}_{m_1, m_2}$  up to $10^9$. 
		}
	\label{fig1}
\end{figure}

We determined the ascending rankings of all ordered pairs of relatively prime coefficients $1\leq m_1 \neq m_2\leq 40$ by $\overline{p^*}_{m_1, m_2}$ up to $10^9$ (Ranking 1), by the maximum of ${p^*}_{m_1, m_2}(n)$ up to $10^9$ (Ranking 2)  and by $R_{m_1, m_2}$ (Ranking 3). 
Subfigures \ref{subfig1}, \ref{subfig2} and \ref{subfig3} in Figure \ref{fig1}  display the  scatter plots of Ranking 3 (vertical axis) against Ranking 1 (horizontal axis), Ranking 2 (vertical axis) against Ranking 1 (horizontal axis), and Ranking 2 (vertical axis) against Ranking 1 (horizontal axis), respectively. Table \ref{TableRank} shows Rankings 1, 2 and 3   of all relatively prime coefficients $1\leq m_1 \neq m_2\leq 20$. The Spearman's rank-order correlation coefficient between $R_{m_1, m_2}$ and $\overline{p^*}_{m_1, m_2}$ up to $10^9$ is  $0.9949$ ($4$ d.p.) and between $R_{m_1, m_2}$ and the maximum of $p^*_{m_1, m_2}(n)$ up to $10^9$ is $0.9958$ ($4$ d.p.). These data suggest very strong rank correlation between $R_{m_1, m_2}$ and $\overline{p^*}_{m_1, m_2}$ up to $10^9$ and also between $R_{m_1, m_2}$ and the maximum of  $p^*_{m_1, m_2}(n)$ up to $10^9$, supporting our hypothesis.

 The function $R_{m_1, m_2}$ also predicts two further observations emerging from our data. Note that the value of $R_{m_1, m_2}$ depends only on $m_2$ and on the set of prime divisors of $m_1$, but not on the actual value of $m_1$. Based on this we can expect $\overline{p^*}_{m_1, m_2}$ and $\overline{p^*}_{m_1', m_2}$ to be `very close' to each other, if $m_1$ and $m_1'$ have the same set of prime divisors. This prediction is in remarkable correspondence with our experimental data, as it is illustrated by Table \ref{Table_m1SamePrFact} in Section \ref{Sect_ToD}. The other observation is the following: For every pair $1\leq m_1< m_2 \leq 40$ relatively prime  we found that $\overline{p^*}_{m_2, m_1} < \overline{p^*}_{m_1, m_2}$ up to $10^9$, which is predicted by the fact that we have $R_{m_2, m_1}<R_{m_1, m_2}$ when $m_1<m_2$.

\begin{figure}[H]
	\centering 
	{\includegraphics{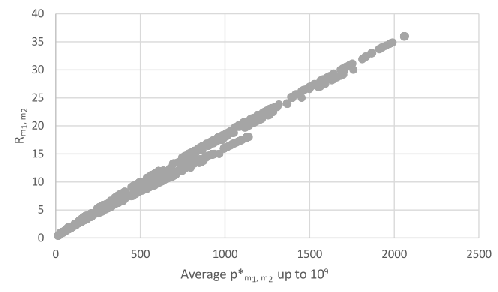}} 
	\caption{$R_{m_1, m_2}$ (vertical axis) plotted against  $\overline{p^*}_{m_1, m_2}$ up to $10^9$ (horizontal axis) for every pair of relatively prime coefficients $1\leq  m_1, m_2 \leq 40$.} \label{Rm1m2avgp}
\end{figure}

Figure \ref{Rm1m2avgp} shows $R_{m_1, m_2}$ plotted against $\overline{p^*}_{m_1, m_2}$ up to $10^9$, for all pairs of relatively prime coefficients $1\leq m_1 \neq m_2\leq 40$. The points are clearly clustered `very close' to a smooth looking curve. The diagram suggests a  relationship close to linear, but we hypothesize that the relationship can be approximated better by some sublinear function, because of the decreasing trend which can be observed in the quotiens $R_{m_1, m_2}/\overline{p^*}_{m_1, m_2}$ as $\overline{p^*}_{m_1, m_2}$ increases. The approximation of this relationship is a topic for further research.

\section{Connections between the magnitude of $p^*_{m_1, m_2}(n)$ and the running times of verifying algorithms of GGC} \label{Sect_Time}

We applied four (two in case $m_1=m_2=1$) different natural algorithms to test GGC  for all pairs of relatively prime coefficients $1\leq m_1, m_2\leq 40$ up to $10^{12}$ \cite{JuhaszBartalos2024}. These algorithms search for $p^*_{m_1, m_2}(n)$ or $p^*_{m_2, m_1}(n)$ in ascending order among the `small' primes or for $q^{**}_{m_1, m_2}(n)$ or $q^{**}_{m_2, m_1}(n)$ in descending order among `large' primes.  According to expectations, for each pair $m_1, m_2$,  algorithms searching for the `large' primes in descending order were faster than those looking for the `small' primes in ascending order.

Based on our data we hypothesize that the significant part of the computation is searching for the $(m_1, m_2)$-partitions of numbers after prime generation.  Therefore, an estimate for the magnitude of the functions $p^*_{m_1, m_2}(n)$ is required to approximate the time complexities of these algorithms.

  When looking for $p^*_{m_1, m_2}(n)$  [$p^*_{m_2, m_1}(n)$], only primes $p$ satisfying 
  \begin{equation}
  m_1p\equiv n \pmod{m_2} \;\;\; [m_2p\equiv n \pmod{m_1}]
  \end{equation} \label{Eqn_Cong}
  were checked as candidates. An analogous idea was also used when searching for $q^{**}_{m_1, m_2}(n)$  [$q^{**}_{m_2, m_1}(n)$].   We introduced the function $$g_{L}(m_1, m_2):= \frac{\overline{p^*}_{L}(m_1, m_2)}{\varphi(m_2)\ln{\overline{p^*}_{L}(m_1, m_2)}},$$        
where $\overline{p^*}_{L}(m_1, m_2)$ is the average of $p^*_{m_1, m_2}(n)$ up to $L$. For $476$  out of all $489$ pairs $m_1<m_2$ tested, we observed that the algorithm searching for $p^*_{m_1, m_2}(n)$ was faster than the one looking for $p^*_{m_2, m_1}(n)$ if and only if 	$g_{10^9}(m_1, m_2)<g_{10^9}(m_2, m_1)$ \cite{JuhaszBartalos2024}.

In less optimized versions of the algorithms where all primes are tested in order as candidates regardless of the congurence conditions \ref{Eqn_Cong}, we predict that if $m_1<m_2$ then searching for  $p^*_{m_2, m_1}(n)$ is faster than for  $p^*_{m_1, m_2}(n)$, since $\overline{p^*}_{L}(m_2, m_1)<\overline{p^*}_{L}(m_1, m_2)$ (where $L=10^9$) in our data, which is also anticipated by $R_{m_2, m_1}<R_{m_1, m_2}$.

\section{Conclusion and closing remarks} \label{Sect_Concl}The function $$ R_{m_1, m_2} := m_2 \cdot \prod_{\substack{s\in P, 2<s\\ s | m_1m_2}}\frac{s-2}{s-1}$$ exhibits very strong rank-order correlation with the average and the maximum of $p^*_{m_1, m_2}(n)$ up to $10^9$ on the pairs $1\leq m_1 \neq m_2\leq 40$ of relatively prime coefficients, with Spearman's rank-order coefficients  $\rho=0.9949$ and $0.9958$,  respectively. By our heuristics we expect this strong rank-order correlation to hold also in general, in cases of all larger thresholds $N\geq 10^9$ and relatively prime positive coefficients $m_1, m_2$. The function $ R_{m_1, m_2}$ helps explain observed trends in our empirical data: (1) For every pair $m_1<m_2$ we found that the average of $p^*_{m_2, m_1}(n)$ up to $10^9$ is less than that of $p^*_{m_1, m_2}(n)$, in line with $R_{m_2, m_1}<R_{m_1, m_2}$.  (2) For every pair $m_1, m_2$ and $m_1', m_2$ such that $m_1$ and $m_1'$ have the same set of prime divisors, the averages of $p^*_{m_1, m_2}(n)$ and $p^*_{m_1', m_2}(n)$, respectively, up to $10^9$ are strikingly close to equal (illustrated by Table \ref{Table_m1SamePrFact}) as predicted by $R_{m_1, m_2}=R_{m_1', m_2}$.  

For the size of the minimal primes in the classical Goldbach-partitions of even numbers Granville \textit{et al.} conjectured $p(n)=O(\log ^2 n\log \log n)$. Granville also proposed two more precise, incompatible conjectures of the form $p(n)\leq (C+o(1))\log^2 n\log{\log{n}}$, where $C$ is `sharp' in the sense that $C$ is the smallest constant with this property: one with $C=C_2^{-1}\approx 1.51478$ and the other one with $C=2e^{-\gamma}C_2^{-1}\approx 1.70098$, where $C_2\approx 0.66016$ is the twin prime constant and $\gamma\approx 0.57722$ is the Euler's constant \cite{Oliveira2014}.
In order to gain information about the running times of the verification algorithms for GGC, estimates of the functions $p^*_{m_1, m_2}(n)$ are required. This is a goal of future research.

\newpage

\section{Tables of data} \label{Sect_ToD}
\nopagebreak

\begin{table}[htbp] 
	\centering
	\caption{The average $\overline{p^*}_{m_1, m_2}$ up to $10^9$ for some pairs $1\leq m_1, m_2 \leq 40$ relatively prime, where $m_2=1, 2\ldots , 10$. The values of $m_1$ in each block of the table have the same set of prime factors.}
	\begin{tabular}{|l|l|l||l|l|l||l|l|l||l|l|l|}
		\hline
		$\mathbf{m_1}$ & $\mathbf{m_2}$ & $\mathbf{\overline{p^*}_{{}_{m_1, m_2}}}$ & $\mathbf{m_1}$ & $\mathbf{m_2}$ & $\mathbf{\overline{p^*}_{{}_{m_1, m_2}}}$ & $\mathbf{m_1}$ & $\mathbf{m_2}$ & $\mathbf{\overline{p^*}_{{}_{m_1, m_2}}}$ & $\mathbf{m_1}$ & $\mathbf{m_2}$ & $\mathbf{\overline{p^*}_{{}_{m_1, m_2}}}$ \\

		&  & up to $10^9$ &  &  & up to $10^9$ &  &  & up to $10^9$ &  &  & up to $10^9$ \\

		\hline
		2     & 1     & 32.80032 & 2     & 3     & 69.3523 & 6     & 5     & 93.49197 & 10    & 7     & 207.4333 \\
		4     & 1     & 32.80632 & 4     & 3     & 69.36289 & 12    & 5     & 93.51101 & 20    & 7     & 207.4348 \\
		8     & 1     & 32.80586 & 8     & 3     & 69.36958 & 18    & 5     & 93.51134 & 40    & 7     & 207.4686 \\  \cline{10-12}
		16    & 1     & 32.80342 & 16    & 3     & 69.36135 & 24    & 5     & 93.49586 & 3     & 8     & 211.8724 \\
		32    & 1     & 32.80091 & 32    & 3     & 69.35155 & 36    & 5     & 93.48327 & 9     & 8     & 211.9202 \\ \cline{1-9}
		3     & 1     & 20.07168 & 5     & 3     & 55.3382 & 14    & 5     & 141.3145 & 27    & 8     & 211.8799 \\           \cline{10-12}
		9     & 1     & 20.07204 & 25    & 3     & 55.32683 & 28    & 5     & 141.2603 & 5     & 8     & 295.392 \\   \cline{4-9}
		27    & 1     & 20.07323 & 10    & 3     & 51.95065 & 5     & 6     & 118.7446 & 25    & 8     & 295.2399 \\  \cline{1-3}  \cline{10-12}
		5     & 1     & 26.76726 & 20    & 3     & 51.9557 & 25    & 6     & 118.7147 & 2     & 9     & 241.7804 \\     \cline{7-9}
		25    & 1     & 26.76609 & 40    & 3     & 51.95532 & 2     & 7     & 277.1068 & 4     & 9     & 241.822 \\   \cline{1-6}
		6     & 1     & 17.27898 & 14    & 3     & 57.04462 & 4     & 7     & 277.1688 & 8     & 9     & 241.7959 \\
		12    & 1     & 17.2793 & 28    & 3     & 57.04528 & 8     & 7     & 277.1186 & 16    & 9     & 241.8264 \\           \cline{4-6}
		18    & 1     & 17.2799 & 3     & 4     & 97.75677 & 16    & 7     & 277.1205 & 32    & 9     & 241.8369 \\        \cline{10-12}
		24    & 1     & 17.28291 & 9     & 4     & 97.77398 & 32    & 7     & 277.0654 & 5     & 9     & 184.588 \\   \cline{7-9}
		36    & 1     & 17.27778 & 27    & 4     & 97.74765 & 3     & 7     & 154.0733 & 25    & 9     & 184.5113 \\   \cline{1-6}     \cline{10-12}
		10    & 1     & 23.95771 & 5     & 4     & 135.3879 & 9     & 7     & 154.1084 & 10    & 9     & 180.7615 \\
		20    & 1     & 23.95715 & 25    & 4     & 135.3506 & 27    & 7     & 154.0745 & 20    & 9     & 180.7731 \\    \cline{4-9}
		40    & 1     & 23.95944 & 2     & 5     & 172.1371 & 5     & 7     & 211.9499 & 40    & 9     & 180.7754 \\   \cline{1-3}   \cline{10-12}
		14    & 1     & 26.55656 & 4     & 5     & 172.1871 & 25    & 7     & 211.871 & 14    & 9     & 198.5837 \\          \cline{7-9}
		28    & 1     & 26.55535 & 8     & 5     & 172.142 & 6     & 7     & 149.3165 & 28    & 9     & 198.5533 \\            \cline{1-3}    \cline{10-12}
		3     & 2     & 43.62564 & 16    & 5     & 172.1673 & 12    & 7     & 149.3576 & 3     & 10    & 208.8865 \\
		9     & 2     & 43.62763 & 32    & 5     & 172.1022 & 18    & 7     & 149.359 & 9     & 10    & 208.9756 \\     \cline{4-6}
		27    & 2     & 43.62843 & 3     & 5     & 97.29165 & 24    & 7     & 149.344 & 27    & 10    & 208.9667 \\            \cline{1-3}
		5     & 2     & 60.48209 & 9     & 5     & 97.31536 & 36    & 7     & 149.3266 &       &       &  \\
		25    & 2     & 60.46878 & 27    & 5     & 97.29562 &       &       &       &       &       &  \\
		\hline
	\end{tabular}
	\label{Table_m1SamePrFact}
\end{table}

\newpage
\pagebreak

\restoregeometry

\tiny{

\begin{landscape}
\hspace*{-5cm} \begin{table}

	\caption{Average (avg) and maximum (max) values of $p^*_{m_1, m_2}(n)$ and $q^*_{m_1, m_2}(n)$ up to $10^9$ for each pair of relatively prime coefficients $1\leq m_1 <  m_2\leq 20$.}\label{tab:b}
	
	\centering

	\begin{tabular}{|l|l|l|l|l|l||l|l|l|l|l|l||l|l|l|l|l|l|} 
		
		\hline
		$\mathbf{m_1}$ & 	$\mathbf{m_2}$ & \multicolumn{2}{c|}{$\mathbf{p^*_{m_1, m_2}(n)}$}  & \multicolumn{2}{c||}{$\mathbf{q^*_{{}_{m_1, m_2}}(n)}$}  & 	$\mathbf{m_1}$ & $\mathbf{m_2}$ & \multicolumn{2}{c|}{$\mathbf{p^*_{m_1, m_2}(n)}$}  & \multicolumn{2}{c||}{$\mathbf{q^*_{{}_{m_1, m_2}}(n)}$}  & 	$\mathbf{m_1}$ & $\mathbf{m_2}$ & \multicolumn{2}{c|}{$\mathbf{p^*_{m_1, m_2}(n)}$}  & \multicolumn{2}{c|}{$\mathbf{q^*_{{}_{m_1, m_2}}(n)}$} \\
		
		 &  & \textbf{avg} & \textbf{max}  & \textbf{avg}  & \textbf{max} &   &   & \textbf{avg}  & \textbf{max}  & \textbf{avg}  & \textbf{max}  &   &   & \textbf{avg}  & \textbf{max}  & \textbf{avg}  & \textbf{max}  \\
		 
		 \hline

		1     & 1     &       &       &       &       & 4     & 9     & 241.822 & 7927  & 97.774 & 3001  & 9     & 13    & 333.584 & 10193 & 222.26 & 6761 \\
		1     & 2     & 80.839 & 3037  & 32.8  & 1609  & 4     & 11    & 494.758 & 19507 & 160.372 & 5939  & 9     & 14    & 331.513 & 10067 & 198.584 & 6337 \\
		1     & 3     & 72.911 & 2371  & 20.072 & 743   & 4     & 13    & 607.515 & 24919 & 163.502 & 6311  & 9     & 16    & 463.174 & 13627 & 241.826 & 7219 \\
		1     & 4     & 181.026 & 6971  & 32.806 & 1453  & 4     & 15    & 327.845 & 9257  & 73.338 & 2153  & 9     & 17    & 464.765 & 13007 & 227.655 & 6481 \\
		1     & 5     & 176.526 & 6833  & 26.767 & 1093  & 4     & 17    & 841.539 & 29669 & 167.531 & 6553  & 9     & 19    & 528.846 & 15649 & 229.533 & 6301 \\
		1     & 6     & 157.484 & 4969  & 17.279 & 643   & 4     & 19    & 960.026 & 32801 & 168.921 & 6947  & 9     & 20    & 449.818 & 13921 & 180.773 & 5519 \\
		1     & 7     & 281.84 & 9431  & 29.376 & 1129  & 5     & 6     & 118.745 & 3457  & 93.492 & 2801  & 10    & 11    & 370.28 & 13093 & 339.42 & 12241 \\
		1     & 8     & 393.604 & 15497 & 32.806 & 1493  & 5     & 7     & 211.95 & 8969  & 145.513 & 6871  & 10    & 13    & 454.483 & 15731 & 345.898 & 11117 \\
		1     & 9     & 245.866 & 8431  & 20.072 & 647   & 5     & 8     & 295.392 & 10369 & 172.142 & 6229  & 10    & 17    & 632.311 & 21647 & 354.305 & 14369 \\
		1     & 10    & 382.522 & 13009 & 23.958 & 1153  & 5     & 9     & 184.588 & 5333  & 97.315 & 2731  & 10    & 19    & 721.599 & 25057 & 357.248 & 13033 \\
		1     & 11    & 500.068 & 17093 & 31.678 & 1499  & 5     & 11    & 375.582 & 11839 & 156.587 & 5881  & 11    & 12    & 304.228 & 8821  & 267.184 & 8293 \\
		1     & 12    & 342.648 & 11261 & 17.279 & 673   & 5     & 12    & 257.838 & 7309  & 93.511 & 2969  & 11    & 13    & 542.8 & 17299 & 452.035 & 16829 \\
		1     & 13    & 612.063 & 23663 & 32.294 & 1297  & 5     & 13    & 459.273 & 16477 & 159.551 & 5521  & 11    & 14    & 542.423 & 20359 & 406.549 & 15227 \\
		1     & 14    & 611.042 & 20359 & 26.557 & 1129  & 5     & 14    & 459.822 & 15773 & 141.315 & 4651  & 11    & 15    & 295.49 & 8941  & 203.796 & 5527 \\
		1     & 15    & 332.373 & 9127  & 15.379 & 557   & 5     & 16    & 638.409 & 24677 & 172.167 & 6451  & 11    & 16    & 754.067 & 26839 & 494.633 & 17863 \\
		1     & 16    & 849.623 & 33997 & 32.803 & 1597  & 5     & 17    & 637.215 & 22751 & 163.446 & 5657  & 11    & 17    & 751.415 & 25621 & 463.029 & 17713 \\
		1     & 17    & 846.422 & 32779 & 33.084 & 1381  & 5     & 18    & 398.036 & 10499 & 93.511 & 2963  & 11    & 18    & 470.391 & 14251 & 267.222 & 7681 \\
		1     & 18    & 529.975 & 15313 & 17.28 & 701   & 5     & 19    & 726.834 & 25609 & 164.784 & 5711  & 11    & 19    & 856.789 & 27581 & 466.872 & 15467 \\
		1     & 19    & 964.977 & 33791 & 33.364 & 1321  & 6     & 7     & 149.317 & 4597  & 129.94 & 3923  & 11    & 20    & 734.055 & 26497 & 370.239 & 12853 \\
		1     & 20    & 825.834 & 29209 & 23.957 & 1069  & 6     & 11    & 267.139 & 8543  & 139.856 & 4813  & 12    & 13    & 328.85 & 9871  & 310.206 & 9479 \\
		2     & 3     & 69.352 & 2083  & 43.626 & 1399  & 6     & 13    & 328.759 & 10883 & 142.486 & 4957  & 12    & 17    & 459.61 & 13033 & 317.598 & 10657 \\
		2     & 5     & 172.137 & 6379  & 60.482 & 2459  & 6     & 17    & 459.546 & 14731 & 145.948 & 4201  & 12    & 19    & 523.692 & 14699 & 320.198 & 9437 \\
		2     & 7     & 277.107 & 12011 & 66.282 & 2663  & 6     & 19    & 523.593 & 16703 & 147.126 & 4423  & 13    & 14    & 552.943 & 19889 & 499.815 & 16843 \\
		2     & 9     & 241.78 & 7129  & 43.628 & 1549  & 7     & 8     & 323.92 & 12589 & 277.119 & 11197 & 13    & 15    & 301.049 & 8539  & 250.574 & 7151 \\
		2     & 11    & 494.633 & 21107 & 71.487 & 3061  & 7     & 9     & 202.501 & 5717  & 154.108 & 4271  & 13    & 16    & 768.659 & 28463 & 607.268 & 25127 \\
		2     & 13    & 607.339 & 21383 & 72.924 & 3049  & 7     & 10    & 315.347 & 9769  & 207.433 & 6841  & 13    & 17    & 765.986 & 25747 & 566.894 & 21851 \\
		2     & 15    & 327.714 & 9049  & 32.917 & 1031  & 7     & 11    & 411.838 & 15131 & 249.973 & 9439  & 13    & 18    & 479.435 & 15199 & 328.849 & 9277 \\
		2     & 17    & 841.438 & 30859 & 74.71 & 3121  & 7     & 12    & 282.761 & 9137  & 149.358 & 4663  & 13    & 19    & 873.509 & 33703 & 571.528 & 22079 \\
		2     & 19    & 959.87 & 34039 & 75.341 & 3001  & 7     & 13    & 504.267 & 18593 & 254.754 & 10099 & 13    & 20    & 748.115 & 27953 & 454.483 & 14851 \\
		3     & 4     & 97.757 & 2939  & 69.363 & 2411  & 7     & 15    & 274.865 & 7499  & 116.406 & 3583  & 14    & 15    & 270.331 & 7789  & 248.994 & 6689 \\
		3     & 5     & 97.292 & 2909  & 55.338 & 1709  & 7     & 16    & 700.594 & 22783 & 277.12 & 10357 & 14    & 17    & 693.436 & 25121 & 566.271 & 20717 \\
		3     & 7     & 154.073 & 4517  & 60.42 & 1789  & 7     & 17    & 698.137 & 24109 & 260.916 & 11069 & 14    & 19    & 791.453 & 27277 & 570.866 & 20873 \\
		3     & 8     & 211.872 & 6869  & 69.37 & 2383  & 7     & 18    & 436.631 & 13367 & 149.359 & 4481  & 15    & 16    & 347.11 & 9521  & 327.84 & 8893 \\
		3     & 10    & 208.887 & 6359  & 51.951 & 1471  & 7     & 19    & 796.433 & 27583 & 263.145 & 10289 & 15    & 17    & 349.899 & 9539  & 308.256 & 8179 \\
		3     & 11    & 271.626 & 8231  & 64.839 & 2113  & 7     & 20    & 682.396 & 23689 & 207.435 & 6841  & 15    & 19    & 398.729 & 10979 & 310.687 & 9109 \\
		3     & 13    & 333.472 & 10733 & 66.064 & 1999  & 8     & 9     & 241.796 & 7027  & 211.92 & 6961  & 16    & 17    & 841.42 & 30727 & 787.381 & 29531 \\
		3     & 14    & 331.38 & 10259 & 57.045 & 1867  & 8     & 11    & 494.594 & 18481 & 348.832 & 13499 & 16    & 19    & 959.865 & 35327 & 793.796 & 28631 \\
		3     & 16    & 463.076 & 13553 & 69.361 & 2239  & 8     & 13    & 607.287 & 23887 & 355.623 & 12107 & 17    & 18    & 491.044 & 14149 & 459.684 & 12953 \\
		3     & 17    & 464.638 & 12503 & 67.628 & 2269  & 8     & 15    & 327.799 & 9091  & 158.333 & 4817  & 17    & 19    & 894.547 & 33721 & 790.894 & 26927 \\
		3     & 19    & 528.697 & 15217 & 68.167 & 2063  & 8     & 17    & 841.438 & 31081 & 364.403 & 15749 & 17    & 20    & 765.917 & 28429 & 632.339 & 25237 \\
		3     & 20    & 449.579 & 12659 & 51.956 & 1579  & 8     & 19    & 959.894 & 42727 & 367.416 & 13999 & 18    & 19    & 523.747 & 14897 & 495.035 & 16943 \\
		4     & 5     & 172.187 & 7109  & 135.388 & 5521  & 9     & 10    & 208.976 & 6469  & 180.762 & 5501  & 19    & 20    & 772.014 & 28729 & 721.715 & 24071 \\
		4     & 7     & 277.169 & 11497 & 148.746 & 5939  & 9     & 11    & 271.709 & 8363  & 218.093 & 6827  &       &       &       &       &       &  \\
		\hline
	\end{tabular}
	\label{Table_AvgMax}
\end{table}
\end{landscape}
}

\newpage

\newpage

\footnotesize{

	\begin{longtable}{|c|c|c|c|c||c|c|c|c|c|}
		\caption{Ascending rankings of all pairs of relatively prime coefficients $1\leq m_1 \neq m_2\leq 20$ by average $p^*_{m_1, m_2}(n)$ and by maximum $p^*_{m_1, m_2}(n)$ up to $10^9$  and by $R_{m_1, m_2}$, respectively.}\label{TableRank} 	\\
		\hline
		$m_1$ & $m_2$ & \multicolumn{3}{c||}{Ranking of $(m_1, m_2)$ by} & $m_1$ & $m_2$ & \multicolumn{3}{c|}{Ranking of $(m_1, m_2)$  by}  \\
		&       &  max $p^*_{m_1, m_2}$  & avg $p^*_{m_1, m_2}$  &  $R_{m_1, m_2}$ &      &       &  max $p^*_{m_1, m_2}$   & avg $p^*_{m_1, m_2}$  & $R_{m_1, m_2}$  \\
		\hline
		\endfirsthead
		
		\hline
		$m_1$ & $m_2$ & \multicolumn{3}{c||}{Ranking of $(m_1, m_2)$ by} & $m_1$ & $m_2$ & \multicolumn{3}{c|}{Ranking of $(m_1, m_2)$  by}  \\
		&       & max $p^*_{m_1, m_2}$  & avg $p^*_{m_1, m_2}$  & $R_{m_1, m_2}$    &      &       & max $p^*_{m_1, m_2}$  & avg $p^*_{m_1, m_2}$  & $R_{m_1, m_2}$   \\
		\hline
		\endhead
		\hline
		\endfoot	
		\hline
		\endlastfoot	
		15    & 1     & 1     & 1     & 1     & 13    & 15    & 115   & 128   & 118 \\
		6     & 1     & 2     & 2     & 4     & 11    & 12    & 117   & 129   & 123 \\
		12    & 1     & 4     & 3     & 4     & 17    & 15    & 110   & 130   & 120 \\
		18    & 1     & 5     & 4     & 4     & 13    & 12    & 131   & 131   & 125 \\
		3     & 1     & 6     & 5     & 4     & 19    & 15    & 123   & 132   & 121 \\
		9     & 1     & 3     & 6     & 4     & 7     & 10    & 134   & 133   & 150 \\
		20    & 1     & 8     & 7     & 8.5   & 17    & 12    & 144   & 134   & 129.5 \\
		10    & 1     & 12    & 8     & 8.5   & 19    & 12    & 129   & 135   & 133 \\
		14    & 1     & 11    & 9     & 11.5  & 7     & 8     & 158   & 136   & 151 \\
		5     & 1     & 9     & 10    & 8.5   & 2     & 15    & 121   & 137   & 129.5 \\
		7     & 1     & 10    & 11    & 11.5  & 8     & 15    & 122   & 138   & 129.5 \\
		11    & 1     & 20    & 12    & 13    & 16    & 15    & 118   & 139   & 129.5 \\
		13    & 1     & 13    & 13    & 14    & 4     & 15    & 126   & 140   & 129.5 \\
		2     & 1     & 24    & 14    & 19.5  & 6     & 13    & 146   & 141   & 143 \\
		16    & 1     & 23    & 15    & 19.5  & 18    & 13    & 127   & 142   & 143 \\
		8     & 1     & 19    & 16    & 19.5  & 12    & 13    & 135   & 143   & 143 \\
		4     & 1     & 17    & 17    & 19.5  & 3     & 14    & 139   & 144   & 137 \\
		15    & 2     & 7     & 18    & 8.5   & 9     & 14    & 136   & 145   & 137 \\
		17    & 1     & 15    & 19    & 15    & 1     & 15    & 124   & 146   & 129.5 \\
		19    & 1     & 14    & 20    & 16    & 3     & 13    & 145   & 147   & 143 \\
		3     & 2     & 16    & 21    & 19.5  & 9     & 13    & 138   & 148   & 143 \\
		9     & 2     & 21    & 22    & 19.5  & 11    & 10    & 156   & 149   & 153.5 \\
		10    & 3     & 18    & 23    & 24    & 1     & 12    & 151   & 150   & 148 \\
		20    & 3     & 22    & 24    & 24    & 13    & 10    & 149   & 151   & 155 \\
		5     & 3     & 25    & 25    & 24    & 15    & 16    & 132   & 152   & 148 \\
		14    & 3     & 27    & 26    & 26.5  & 11    & 8     & 168   & 153   & 158 \\
		7     & 3     & 26    & 27    & 26.5  & 15    & 17    & 133   & 154   & 146 \\
		5     & 2     & 38    & 28    & 35    & 17    & 10    & 175   & 155   & 156 \\
		11    & 3     & 31    & 29    & 28    & 13    & 8     & 155   & 156   & 159 \\
		13    & 3     & 28    & 30    & 29    & 19    & 10    & 165   & 157   & 157 \\
		7     & 2     & 39    & 31    & 39    & 17    & 8     & 189   & 158   & 165 \\
		17    & 3     & 34    & 32    & 30    & 19    & 8     & 172   & 159   & 168 \\
		19    & 3     & 29    & 33    & 31    & 20    & 11    & 160   & 160   & 161 \\
		2     & 3     & 30    & 34    & 35    & 10    & 11    & 166   & 161   & 161 \\
		16    & 3     & 33    & 35    & 35    & 5     & 11    & 153   & 162   & 161 \\
		4     & 3     & 37    & 36    & 35    & 1     & 10    & 163   & 163   & 165 \\
		8     & 3     & 36    & 37    & 35    & 1     & 8     & 186   & 164   & 175 \\
		11    & 2     & 50    & 38    & 40    & 5     & 18    & 143   & 165   & 153.5 \\
		1     & 3     & 35    & 39    & 35    & 15    & 19    & 147   & 166   & 152 \\
		13    & 2     & 49    & 40    & 41    & 14    & 11    & 183   & 167   & 179 \\
		15    & 4     & 32    & 41    & 35    & 7     & 11    & 180   & 168   & 179 \\
		17    & 2     & 51    & 42    & 44.5  & 7     & 18    & 167   & 169   & 165 \\
		19    & 2     & 47    & 43    & 48    & 3     & 20    & 159   & 170   & 165 \\
		1     & 2     & 48    & 44    & 50    & 9     & 20    & 171   & 171   & 165 \\
		6     & 5     & 41    & 45    & 44.5  & 13    & 11    & 193   & 172   & 193 \\
		12    & 5     & 45    & 46    & 44.5  & 10    & 13    & 188   & 173   & 185 \\
		18    & 5     & 44    & 47    & 44.5  & 20    & 13    & 178   & 174   & 185 \\
		3     & 5     & 42    & 48    & 44.5  & 5     & 13    & 191   & 175   & 185 \\
		9     & 5     & 40    & 49    & 44.5  & 6     & 17    & 177   & 176   & 171 \\
		3     & 4     & 43    & 50    & 50    & 12    & 17    & 164   & 177   & 171 \\
		9     & 4     & 46    & 51    & 50    & 18    & 17    & 161   & 178   & 171 \\
		15    & 7     & 53    & 52    & 52    & 5     & 14    & 190   & 179   & 183 \\
		5     & 6     & 52    & 53    & 53    & 17    & 11    & 198   & 180   & 194 \\
		7     & 6     & 54    & 54    & 54    & 3     & 16    & 169   & 181   & 175 \\
		5     & 4     & 70    & 55    & 65    & 9     & 16    & 170   & 182   & 175 \\
		11    & 6     & 63    & 56    & 55    & 3     & 17    & 157   & 183   & 171 \\
		14    & 5     & 61    & 57    & 67.5  & 9     & 17    & 162   & 184   & 171 \\
		13    & 6     & 65    & 58    & 56    & 19    & 11    & 185   & 185   & 195 \\
		7     & 5     & 96    & 59    & 67.5  & 11    & 18    & 174   & 186   & 177 \\
		17    & 6     & 55    & 60    & 57    & 13    & 18    & 181   & 187   & 179 \\
		19    & 6     & 57    & 61    & 58    & 17    & 18    & 173   & 188   & 181 \\
		7     & 4     & 77    & 62    & 69    & 8     & 11    & 200   & 189   & 198 \\
		6     & 7     & 60    & 63    & 61    & 2     & 11    & 208   & 190   & 198 \\
		12    & 7     & 62    & 64    & 61    & 16    & 11    & 199   & 191   & 198 \\
		18    & 7     & 58    & 65    & 61    & 4     & 11    & 202   & 192   & 198 \\
		3     & 7     & 59    & 66    & 61    & 19    & 18    & 195   & 193   & 182 \\
		9     & 7     & 56    & 67    & 61    & 14    & 13    & 194   & 194   & 201.5 \\
		11    & 5     & 76    & 68    & 71.5  & 1     & 11    & 196   & 195   & 198 \\
		1     & 6     & 66    & 69    & 65    & 7     & 13    & 201   & 196   & 201.5 \\
		15    & 8     & 64    & 70    & 65    & 6     & 19    & 192   & 197   & 189 \\
		13    & 5     & 71    & 71    & 74    & 12    & 19    & 176   & 198   & 189 \\
		11    & 4     & 78    & 72    & 77    & 18    & 19    & 179   & 199   & 189 \\
		17    & 5     & 73    & 73    & 75    & 3     & 19    & 182   & 200   & 189 \\
		13    & 4     & 81    & 74    & 78    & 9     & 19    & 187   & 201   & 189 \\
		19    & 5     & 74    & 75    & 76    & 1     & 18    & 184   & 202   & 192 \\
		17    & 4     & 88    & 76    & 84.5  & 11    & 14    & 205   & 203   & 203 \\
		19    & 4     & 97    & 77    & 90    & 11    & 13    & 197   & 204   & 205 \\
		2     & 5     & 84    & 78    & 84.5  & 13    & 14    & 203   & 205   & 204 \\
		8     & 5     & 79    & 79    & 84.5  & 17    & 14    & 206   & 206   & 206 \\
		16    & 5     & 85    & 80    & 84.5  & 17    & 13    & 211   & 207   & 208 \\
		4     & 5     & 101   & 81    & 84.5  & 19    & 14    & 207   & 208   & 207 \\
		1     & 5     & 92    & 82    & 84.5  & 19    & 13    & 212   & 209   & 209 \\
		10    & 9     & 68    & 83    & 71.5  & 16    & 13    & 224   & 210   & 213 \\
		20    & 9     & 69    & 84    & 71.5  & 8     & 13    & 217   & 211   & 213 \\
		1     & 4     & 99    & 85    & 92    & 2     & 13    & 209   & 212   & 213 \\
		5     & 9     & 67    & 86    & 71.5  & 4     & 13    & 221   & 213   & 213 \\
		14    & 9     & 82    & 87    & 84.5  & 1     & 14    & 204   & 214   & 210 \\
		7     & 9     & 75    & 88    & 84.5  & 1     & 13    & 215   & 215   & 213 \\
		15    & 11    & 72    & 89    & 79    & 10    & 17    & 210   & 216   & 217 \\
		10    & 7     & 93    & 90    & 99.5  & 20    & 17    & 225   & 217   & 217 \\
		20    & 7     & 94    & 91    & 99.5  & 5     & 17    & 213   & 218   & 217 \\
		3     & 10    & 83    & 92    & 84.5  & 5     & 16    & 220   & 219   & 219 \\
		9     & 10    & 86    & 93    & 84.5  & 7     & 20    & 216   & 220   & 220 \\
		3     & 8     & 95    & 94    & 92    & 14    & 17    & 223   & 221   & 221.5 \\
		9     & 8     & 98    & 95    & 92    & 7     & 17    & 219   & 222   & 221.5 \\
		5     & 7     & 120   & 96    & 99.5  & 7     & 16    & 214   & 223   & 223 \\
		11    & 9     & 91    & 97    & 94    & 10    & 19    & 222   & 224   & 225 \\
		13    & 9     & 90    & 98    & 95    & 20    & 19    & 218   & 225   & 225 \\
		17    & 9     & 87    & 99    & 96    & 5     & 19    & 226   & 226   & 225 \\
		19    & 9     & 80    & 100   & 97    & 11    & 20    & 229   & 227   & 227 \\
		2     & 9     & 102   & 101   & 105.5 & 13    & 20    & 235   & 228   & 228 \\
		8     & 9     & 100   & 102   & 105.5 & 11    & 17    & 227   & 229   & 231 \\
		4     & 9     & 109   & 103   & 105.5 & 11    & 16    & 230   & 230   & 232 \\
		16    & 9     & 104   & 104   & 105.5 & 17    & 20    & 236   & 231   & 229 \\
		1     & 9     & 114   & 105   & 105.5 & 13    & 17    & 228   & 232   & 233 \\
		15    & 14    & 89    & 106   & 99.5  & 13    & 16    & 237   & 233   & 234 \\
		11    & 7     & 130   & 107   & 119   & 19    & 20    & 239   & 234   & 230 \\
		15    & 13    & 103   & 108   & 102   & 17    & 16    & 241   & 235   & 237.5 \\
		13    & 7     & 137   & 109   & 122   & 19    & 17    & 231   & 236   & 239 \\
		5     & 12    & 105   & 110   & 105.5 & 14    & 19    & 232   & 237   & 235.5 \\
		17    & 7     & 148   & 111   & 124   & 19    & 16    & 238   & 238   & 240 \\
		19    & 7     & 140   & 112   & 126   & 7     & 19    & 234   & 239   & 235.5 \\
		6     & 11    & 116   & 113   & 113   & 1     & 20    & 240   & 240   & 237.5 \\
		12    & 11    & 112   & 114   & 113   & 16    & 17    & 243   & 241   & 243 \\
		18    & 11    & 107   & 115   & 113   & 2     & 17    & 244   & 242   & 243 \\
		14    & 15    & 108   & 116   & 109.5 & 8     & 17    & 245   & 243   & 243 \\
		3     & 11    & 111   & 117   & 113   & 4     & 17    & 242   & 244   & 243 \\
		9     & 11    & 113   & 118   & 113   & 1     & 17    & 246   & 245   & 243 \\
		7     & 15    & 106   & 119   & 109.5 & 1     & 16    & 251   & 246   & 246 \\
		2     & 7     & 154   & 120   & 137   & 11    & 19    & 233   & 247   & 247 \\
		8     & 7     & 150   & 121   & 137   & 13    & 19    & 248   & 248   & 248 \\
		16    & 7     & 141   & 122   & 137   & 17    & 19    & 249   & 249   & 249 \\
		4     & 7     & 152   & 123   & 137   & 16    & 19    & 253   & 250   & 252 \\
		1     & 7     & 128   & 124   & 137   & 2     & 19    & 252   & 251   & 252 \\
		7     & 12    & 125   & 125   & 116   & 8     & 19    & 254   & 252   & 252 \\
		5     & 8     & 142   & 126   & 148   & 4     & 19    & 247   & 253   & 252 \\
		11    & 15    & 119   & 127   & 117   & 1     & 19    & 250   & 254   & 252 \\
\end{longtable}
}

\newpage

\bibliographystyle{amsplain}

\end{document}